# Ratio tests for change point detection

**Lajos Horváth**[*1], **Zsuzsanna Horváth**[1] **and Marie Hušková**[†2]

*University of Utah and Charles University*

**Abstract:** We propose new tests to detect a change in the mean of a time series. Like many existing tests, the new ones are based on the CUSUM process. Existing CUSUM tests require an estimator of a scale parameter to make them asymptotically distribution free under the no change null hypothesis. Even if the observations are independent, the estimation of the scale parameter is not simple since the estimator for the scale parameter should be at least consistent under the null as well as under the alternative. The situation is much more complicated in case of dependent data, where the empirical spectral density at 0 is used to scale the CUSUM process. To circumvent these difficulties, new tests are proposed which are ratios of CUSUM functionals. We demonstrate the applicability of our method to detect a change in the mean when the errors are AR(1) and GARCH(1,1) sequences.

## 1. Introduction

Change point detection is an important part of statistical and economic analysis. Predictions and statistical inference will be invalid if changes in the regimes during the data collection period are not taken into account. The main problems in the change point analysis are to decide whether the statistical model for a series of observations does not change (no change situation) or whether the model changes one or more times and in the latter case to identify when the changes have occured. For surveys on change point methods we refer to Csörgő and Horváth [3] and Perron [11].

In this paper we consider at most one change in the location model

$$X_k = \mu_k + \epsilon_k \qquad 1 \le k \le n,$$

where $\mu_1, \dots, \mu_n$ are the means of the respective observations while $\epsilon_1, \dots, \epsilon_n$ are random error terms with zero mean satisfying some additional assumptions specified below. Under the no change null hypothesis

$$H_0: \quad \mu_k = \mu \quad 1 \le k \le n$$

while under the alternative

$$H_A: \text{ there is } 1 \le k^* < n \text{ such that } \mu_1 = \mu_2 = \cdots = \mu_k \ne \mu_{k+1} = \cdots = \mu_n.$$

*Supported by NSF Grant DMS-06-04670 and Grant RGC-HKUST6428/06H.
†Supported in part by Grants MSM 02162839, GAČR 201/06/0186 and LC 06024.
[1]Department of Mathematics, University of Utah, 155 South 1440 East, Salt Lake City, UT 84112-0090, USA, e-mail: Horvath@math.utah.edu; c-hhz@math.utah.edu
[2]Department of Statistics, Charles University, Sokolovská 83, 18675 Prague, Czech Republic, e-mail: huskova@karlin.mff.cuni.cz

The most popular methods are based on functionals of properly standardized cumulative sums (CUSUM) $\sum_{i=1}^{k}(X_i - \bar{X}_n)$, $k = 1, \ldots, n$, where $\bar{X}_n = (1/n)\sum_{1\le i\le n} X_i$. For example, we reject $H_0$ if

$$T_{n,1} = \max_{1\le k\le n} |\sum_{i=1}^{k}(X_i - \bar{X}_n)|/(n\sigma_n^2)^{1/2}$$

is large, where and $n\sigma_n^2$ is a proper estimator for the variance of $\sum_{1\le i\le n} \epsilon_i$. Similarly, the $R/S$ statistic proposed by Lo [10] is

$$T_{n,2} = \frac{1}{(n\sigma_n^2)^{1/2}}\left[\max_{1\le k\le n}\sum_{i=1}^{k}(X_i - \bar{X}_n) - \min_{1\le k\le n}\sum_{i=1}^{k}(X_i - \bar{X}_n)\right].$$

The statistic $T_{n,2}$ was modified by Giraitis et al. [5] who introduced

$$T_{n,3} = \frac{1}{n^2\sigma_n^2}\left(\sum_{k=1}^{n}\left[\sum_{j=1}^{k}(X_j - \bar{X}_n)\right]^2 - \frac{1}{n}\left[\sum_{k=1}^{n}\sum_{j=1}^{k}(X_j - \bar{X}_n)\right]^2\right).$$

Under suitable assumptions on the error terms all three test statistics are sensitive w.r.t. change(s) in the mean (location). Asymptotic properties of $T_{n,1}, i = 1, 2, 3$ were derived under the conditions

$$E\epsilon_i = 0, \quad 1 \le i < \infty \tag{1.1}$$

and

$$\text{there is } \sigma > 0 \text{ such that } n^{-1/2}\sum_{1\le i\le nt}\epsilon_i \xrightarrow{\mathcal{D}[0,1]} \sigma W(t) \tag{1.2}$$

where $\{W(t), 0 \le t < \infty\}$ is a Wiener process. ($\xrightarrow{\mathcal{D}[0,1]}$ denotes weak convergence in $\mathcal{D}[0,1]$.) Condition (1.2) means that $\{\epsilon_i\}$ is a weakly dependent sequence which satisfies the functional central limit theorem. Since $T_{n,i}, i = 1, 2, 3$ are functions of $n^{-1/2}\sum_{1\le i\le nt}\epsilon_i$, (1.2) will yield the asymptotic distributions of the test statistics both under $H_0$ and $H_A$. However, the estimator $\sigma_n^2$ must satisfy that $\sigma_n^2 \xrightarrow{P} \sigma^2$ under $H_0$ and at least it must be bounded in probability under the alternative. If $\{\epsilon_i\}$ is a strictly stationary sequence with $0 < E\epsilon_0^2 < \infty$, then the estimation of $\sigma^2$ is based on the fact that it is related to the spectral density at 0. So one needs to choose a kernel and the number of lags used in the estimation. One of the most popular choices is Bartlett's estimator. However the rate of convergence is very slow even under $H_0$ and $\sigma_n^2$ might go to infinity under $H_A$. Modifications of the Bartlett estimator in the change point context can be found in Berkes et al. [1] and Berkes et al. [2]. Therefore it is desirable to develop procedures for testing $H_0$ against $H_A$, where the estimator of $\sigma^2$ from (1.2) is not needed. We develop such test procedures based on functionals of CUSUMs.

In the definitions of $T_{n,i}, i = 1, 2, 3$, functionals of CUSUMs are computed for the first $k$ and the last $n - k$ observations. If the difference between functionals is large for at least one $k$, the null hypothesis of no change is rejected. We suggest computing the ratio of the CUSUM functionals instead of the differences. This way there will be no need for the estimation of $\sigma^2$. Instead of using $T_{n,1}$ we suggest

$$V_{n,1} = \max_{n\delta\le k\le n-n\delta}\frac{\max_{1\le i\le k}\left|\sum_{1\le j\le i}(X_j - \bar{X}_k)\right|}{\max_{k\le i\le n}\left|\sum_{i\le j\le n}(X_j - \tilde{X}_k)\right|},$$

where $0 < \delta < 1/2$ and

$$\bar{X}_k = \frac{1}{k} \sum_{1 \le i \le k} X_i \qquad \text{and} \quad \tilde{X}_k = \frac{1}{n-k} \sum_{k < i \le n} X_i.$$

Similarly,

$$V_{n,2} = \max_{n\delta \le k \le n - n\delta} \frac{\max_{1\le i\le k} \sum_{1\le j\le i}(X_j - \bar{X}_k) - \min_{1\le i\le k} \sum_{1\le j\le i}(X_j - \bar{X}_k)}{\max_{k<i\le n} \sum_{i\le j\le n}(X_j - \tilde{X}_k) - \min_{k<i\le n} \sum_{i\le j\le n}(X_j - \tilde{X}_k)}$$

and

$$V_{n,3} = \max_{n\delta \le k \le n-\delta n} \frac{\sum_{i=1}^k \left[\sum_{j=1}^i (X_j - \bar{X}_k)\right]^2 - \frac{1}{k}\left[\sum_{i=1}^k \sum_{j=1}^i (X_j - \bar{X}_k)\right]^2}{\sum_{i=k+1}^n \left[\sum_{j=i}^n (X_j - \tilde{X}_k)\right]^2 - \frac{1}{n-k}\left[\sum_{i=k+1}^n \sum_{j=i}^n (X_j - \tilde{X}_k)\right]^2}.$$

Our first result gives the convergence in distribution results for $V_{n,i}, i = 1, 2, 3$ under $H_0$. Let $W(t), 0 \le t < \infty$ be a Wiener process and define the following processes:

$$\eta_{1,1}(t) = \sup_{0\le s\le t} |W(s) - (s/t)W(t)|,$$

$$\eta_{1,2}(t) = \sup_{t\le s\le 1} |W^*(s) - (1-s)/(1-t)W^*(t)|,$$

$$\eta_{2,1}(t) = \sup_{0\le s\le t} (W(s) - (s/t)W(t)) - \inf_{0\le s\le t}(W(s) - (s/t)W(t)),$$

$$\eta_{2,2}(t) = \sup_{t\le s\le 1} (W^*(s) - ((1-s)/(1-t))W^*(t))$$

$$- \inf_{t\le s\le 1} (W^*(s) - ((1-s)/(1-t))W^*(t)),$$

$$\eta_{3,1}(t) = \int_0^t (W(s) - (s/t)W(t))^2 ds - \frac{1}{t}\left(\int_0^t (W(s) - (s/t)W(t))ds\right)^2$$

and

$$\eta_{3,2}(t) = \int_t^1 (W^*(s) - ((1-s)/(1-t))W^*(t))^2 ds$$

$$-\frac{1}{1-t}\left(\int_t^1 (W^*(s) - ((1-s)/(1-t))W^*(t))ds\right)^2,$$

where $W^*(t) = W(1) - W(t)$.

**Theorem 1.1.** *If $H_0$, (1.1) and (1.2) hold then*

$$V_{n,1} \xrightarrow{\mathcal{D}} \sup_{\delta\le t\le 1-\delta} \frac{\eta_{1,1}(t)}{\eta_{1,2}(t)}, \tag{1.3}$$

$$V_{n,2} \xrightarrow{\mathcal{D}} \sup_{\delta\le t\le 1-\delta} \frac{\eta_{2,1}(t)}{\eta_{2,2}(t)} \tag{1.4}$$

*and*

$$V_{n,3} \xrightarrow{\mathcal{D}} \sup_{\delta\le t\le 1-\delta} \frac{\eta_{3,1}(t)}{\eta_{3,2}(t)}. \tag{1.5}$$

We note that $\sup_{\delta\le t\le 1-\delta}$ can be replaced with $\sup_{0<t\le 1-\delta}$ in (1.3)–(1.5). Since

$$\lim_{t\to 1-}\sup_{t\le s\le 1}\left|W^*(s)-\frac{1-s}{1-t}W^*(t)\right|=0 \ \text{ a.s.},$$

$$\lim_{t\to 1-}\sup_{t\le s\le 1}\left(W^*(s)-\frac{1-s}{1-t}W^*(t)\right) - \inf_{t\le s\le 1}\left(W^*(s)-\frac{1-s}{1-t}W^*(t)\right)=0 \ \text{ a.s.,} \tag{1.6}$$

and

$$\lim_{t\to 1-}\int_t^1\left(W^*(s)-\frac{1-s}{1-t}W^*(t)\right)^2 ds - \frac{1}{1-t}\left(\int_t^1\left(W^*(s)-\frac{1-s}{1-t}W^*(t)\right)ds\right)^2=0 \quad \text{a.s.}$$

we cannot replace $\sup_{\delta\le t\le 1-\delta}$ with $\sup_{\delta\le t<1}$, in (1.3)–(1.5).

The Wiener process $W$ has independent increments and therefore for any $0<t<1$ we have that $\{W(s)-(s/t)W(t), 0\le s\le t\}$ and $\{W^*(s)-((1-s)/(1-t))W^*(t), t\le s\le 1\}$ are independent. Change of variable and the scale transformation of $W$ give that

$$\eta_{1,1}(t)=\sup_{0\le u\le 1}|W(ut)-uW(t)| \overset{\mathcal{D}}{=} t^{1/2}\sup_{0\le u\le 1}|B(u)|, \quad \text{for all} \quad 0<t<1,$$

where $B(u)=W(u)-uW(1)$ is a Brownian bridge. Therefore for any $0<t<1$

$$\frac{\eta_{1,1}(t)}{\eta_{1,2}(t)} \overset{\mathcal{D}}{=} \left(\frac{t}{1-t}\right)^{1/2}\frac{\sup_{0\le u\le 1}|B_1(u)|}{\sup_{0\le u\le 1}|B_2(u)|},$$

where $\{B_1(u), 0\le u\le 1\}$ and $\{B_2(u), 0\le u\le 1\}$ are independent Brownian bridges. Similar arguments give

$$\frac{\eta_{2,1}(t)}{\eta_{2,2}(t)} \overset{\mathcal{D}}{=} \left(\frac{t}{1-t}\right)^{1/2}\frac{\sup_{0\le u\le 1}B_1(u)-\inf_{0\le u\le 1}B_1(u)}{\sup_{0\le u\le 1}B_2(u)-\inf_{0\le u\le 1}B_2(u)}$$

and

$$\frac{\eta_{3,1}(t)}{\eta_{3,2}(t)} \overset{\mathcal{D}}{=} \frac{t}{1-t}\,\frac{\int_0^1 B_1^2(u)du-\left(\int_0^1 B_1(u)du\right)^2}{\int_0^1 B_2^2(u)du-\left(\int_0^1 B_2(u)du\right)^2}$$

for any $0<t<1$.

Kim [7] used ratio tests to detect changes in the persistence of a linear time series. The asymptotic as well as the finite sample properties (including power) of Kim's test were investigated by Kim et al. [8] and Leybourne and Taylor [9]. Since we try to detect changes in the means (location) our tests are different from Kim's so we need to investigate the asymptotic power. Let $\Delta_n=\mu_{k^*}-\mu_{k^*+1}$ be the size of the change.

**Theorem 1.2.** *If (1.1) and (1.2) and $H_A$ hold,*

$$k^*=[n\theta] \ \textit{ with some } \ 0<\theta<1, \tag{1.7}$$

*and*

$$n^{1/2}|\Delta_n| \to \infty, \tag{1.8}$$

*then*

$$V_{n,i} \xrightarrow{P} \infty, \quad i = 1, 2, 3,$$

*assuming that*

$$\delta < \theta < 1 - \delta. \tag{1.9}$$

Our tests were developed to check if the mean has changed at an unknown time. Ratio type tests can also be used to see if the sequence changes from "stationary" into "difference stationary" . We say that the sequence is "stationary" if the sum of the $X_k$'s satisfies the functional central limit theorem and "difference stationary" if the $X_k$'s themselves satisfy the functional central limit theorem with suitable normalization. Now we consider the following alternative:

$$H_A^*: \quad \text{there is } 1 \le k^* < n \text{ such that } X_k = \begin{cases} \mu + \epsilon_k, & 1 \le k \le k^*, \\ \mu + \epsilon_{k^*} + \ldots + \epsilon_k, & k^* < k \le n. \end{cases}$$

The statistics $V_{n,i}, i = 1, 2, 3$ may not be able to detect if $H_0$ or $H_A^*$ hold, since even under $H_A^*$ the statistics have nondegenerate limit distributions. Hence we suggest the following modification of the test statistics to detect $H_A^*$:

$$Z_{n,1} = \max_{n\delta \le k \le n - n\delta} \frac{\max_{k \le i \le n} \left| \sum_{i \le j \le n} (X_j - \tilde{X}_k) \right|}{\max_{1 \le i \le k} \left| \sum_{1 \le j \le i} (X_j - \bar{X}_k) \right|},$$

$$Z_{n,2} = \max_{n\delta \le k \le n - n\delta} \frac{\max_{k < i \le n} \sum_{i \le j \le n} (X_j - \tilde{X}_k) - \min_{k < i \le n} \sum_{i \le j \le n} (X_j - \tilde{X}_k)}{\max_{1 \le i \le k} \sum_{1 \le j \le i} (X_j - \bar{X}_k) - \min_{1 \le i \le k} \sum_{1 \le j \le i} (X_j - \bar{X}_k)}$$

and

$$Z_{n,3} = \max_{n\delta \le k \le n - \delta n} \frac{\sum_{i=k+1}^{n} \left[ \sum_{j=i}^{n} (X_j - \tilde{X}_k) \right]^2 - \frac{1}{n-k} \left[ \sum_{i=k+1}^{n} \sum_{j=i}^{n} (X_j - \tilde{X}_j) \right]^2}{\sum_{i=1}^{k} \left[ \sum_{j=1}^{i} (X_j - \bar{X}_k) \right]^2 - \frac{1}{k} \left[ \sum_{i=1}^{k} \sum_{j=1}^{i} (x_j - \bar{X}_k) \right]^2}$$

The limit distributions of $Z_{n,i}, i = 1, 2, 3$ can be easily derived following the proof of Theorem 1.2 and we get the following results:

$$Z_{n,1} \xrightarrow{\mathcal{D}} \sup_{\delta \le t \le 1-\delta} \frac{\eta_{1,2}(t)}{\eta_{1,1}(t)}, \quad Z_{n,2} \xrightarrow{\mathcal{D}} \sup_{\delta \le t \le 1-\delta} \frac{\eta_{2,2}(t)}{\eta_{2,1}(t)} \quad \text{and} \quad Z_{n,3} \xrightarrow{\mathcal{D}} \sup_{\delta \le t \le 1-\delta} \frac{\eta_{3,2}(t)}{\eta_{3,1}(t)}.$$

**Theorem 1.3.** *If (1.1), (1.2), (1.7)–(1.9) and $H_A^*$ hold, then*

$$Z_{n,i} \xrightarrow{P} \infty, \quad i = 1, 2, 3.$$

However, $V_{n,i}, i = 1, 2, 3$ have power against the alternative

$$H_A^{**}: \quad \text{there is } 1 \le k^* < n \text{ such that } X_k = \begin{cases} \mu + \epsilon_{k^*} + \cdots + \epsilon_k, & 1 \le k \le k^*, \\ \mu + \epsilon_k, & k^* < k \le n. \end{cases}$$

**Theorem 1.4.** *If (1.1), (1.2), (1.7)–(1.9) and $H_A^{**}$ hold, then*

$$V_{n,i} \xrightarrow{P} \infty, \quad i = 1, 2, 3.$$

The alternative $H_A^{**}$ is somewhat the opposite of $H_A^*$; the first observations follow a random walk and at $k^*$ they turn into a "stationary" sequence. Of course, the statistics may not detect the difference between $H_0$ and $H_A^*$. If we are interested only if a change occured from or into a random walk at an unknown time, i/e.,we are testing $H_0$ against $H_A^* \cup H_A^{**}$, we must combine $V_{n,i}$ and $Z_{n,i}$. Let

$$\widetilde{T}_{n,i} = \max(V_{n,i}, Z_{n,i}) \;\; i = 1, 2, 3.$$

Following the proof of Theorem 1.1 one can easily verify that under $H_0$

$$\widetilde{T}_{n,1} \xrightarrow{\mathcal{D}} \max\left\{ \sup_{\delta \le t \le 1-\delta} \frac{\eta_{1,2}(t)}{\eta_{1,1}(t)}, \sup_{\delta \le t \le 1-\delta} \frac{\eta_{1,1}(t)}{\eta_{1,2}(t)} \right\},$$
$$\widetilde{T}_{n,2} \xrightarrow{\mathcal{D}} \max\left\{ \sup_{\delta \le t \le 1-\delta} \frac{\eta_{2,2}(t)}{\eta_{2,1}(t)}, \sup_{\delta \le t \le 1-\delta} \frac{\eta_{2,1}(t)}{\eta_{2,2}(t)} \right\},$$
$$T_{n,3} \xrightarrow{\mathcal{D}} \max\left\{ \sup_{\delta \le t \le 1-\delta} \frac{\eta_{3,2}(t)}{\eta_{3,1}(t)}, \sup_{\delta \le t \le 1-\delta} \frac{\eta_{3,1}(t)}{\eta_{3,2}(t)}, \right\}.$$

**Theorem 1.5.** *If (1.1), (1.2), (1.7)–(1.9) and $H_A^{**}$ hold, then*

$$\widetilde{T}_{n,i} \xrightarrow{P} \infty, \quad i = 1, 2, 3.$$

**Remark 1.1.** It would be more natural to use $\delta = 0$ in our results. However, as we pointed out after Theorem 1.1, (1.6) yields that $\sup_{0<t<1} \eta_{1,1}(t)/\eta_{1,2}(t) = \infty$ a.s. By Chung's law (cf. Csörgő and Révész [4]) for any $\nu > 1/2$

$$\lim_{t \to 1-} (1-t)^{-\nu} \eta_{1,2}(t) = \infty \quad \text{a.s.}$$

and therefore

$$P\{ \sup_{0<t<1} (1-t)^{\nu} \eta_{1,1}(t)/\eta_{1,2}(t) < \infty \} = 1. \tag{1.10}$$

By (1.10) we conjecture that

$$\max_{1<k<n} \left(1 - \frac{k}{n}\right)^{\nu} \frac{\max_{1 \le i \le k} | \sum_{1 \le j \le i} (X_j - \bar{X}_k)|}{\max_{k \le i \le n} | \sum_{i \le j \le n} (X_j \tilde{X}_k)|} \xrightarrow{\mathcal{D}} \sup_{0<t<1} (1-t)^{\nu} \frac{\eta_{1,1}(t)}{\eta_{1,2}(t)}.$$

Using weight functions, the statistics $V_{n,2}$ and $V_{n,3}$ can be modified in a similar way so one can take $\delta = 0$ in the weighted statistics.

**Remark 1.2.** We would like to note that ratio tests can be derived not only for partial sums with a Wiener limit but for more general processes.

## 2. Applications

In our first example the error terms are linear processes defined as

$$\epsilon_k = \sum_{i=0}^{\infty} \alpha_i \delta_{k-i}, \tag{2.1}$$

where
(2.2)
$\delta_i, -\infty < i < \infty$ are independent identically distributed random variables

and

$$E\delta_0 = 0 \text{ and } E\delta_0^2 < \infty.$$

If

$$\sum_{0\le i<\infty} |\alpha_i| < \infty \text{ and } \sum_{0\le i<\infty} \alpha_i \ne 0, \tag{2.3}$$

then (1.2) holds with $\sigma^2 = E\delta_0^2(\sum_{0\le i<\infty} \alpha_i)^2$. The proof of (2.3) is in Hannan [6] (cf. also Wang et al. [12]). We would like to note that by Wu and Min [13], (1.2) holds for sums of linear processes without assuming (2.2).

We studied the behaviour of $V_{n,1}$ when $\delta = 0.2$. We used Monte Carlo simulations to get critical values for $\sup_{.2\le t\le .8} \eta_{1,1}(t)/\eta_{1,2}(t)$. In our simulation study we assumed that $\delta_i, -\infty < i < \infty$ are independent standard normal random variables and $c_i = \rho^i$. This means that the observations are elements of of a stationary AR(1) process with parameter $\rho$. The results in Table 1 suggest that the asymptotic critical values are acceptable even for moderate sample sizes, if $\rho$ is not close to 1. The power in Tables 2–7 is a decreasing function of $\rho$ as $\rho$ tends to 1. The location of the time of the change has little effect on the power; the power is nearly the same for $k^* = n/2$ and $k^* = n/4$.

In the second example we assume that $\epsilon_k$ are elements of a GARCH(1,1) sequence. This means that $\epsilon_k$ satisfies the recursion

$$\epsilon_k = \delta_k \tau_k \quad \text{and} \quad \tau_k^2 = \omega + \alpha\epsilon_{k-1}^2 + \beta\tau_{k-1}^2, \quad -\infty < k < \infty,$$

where $\omega > 0, \alpha \ge 0, \beta \ge 0$. Assuming that (2.2) holds and

$$E\delta_0^2 < \infty \text{ and } \alpha E\delta_0^2 + \beta < 1,$$

then Berkes et al. [2] proved that (1.1) is satisfied with $\sigma^2 = \omega/(1 - \alpha E\delta_0^2 - \beta)$.

TABLE 1
*Simulated significance levels for $V_{n,1}$ when $\delta = 0.2$*

| **n** | **Level** | $\rho$: **0.1** | **0.2** | **0.3** | **0.4** | **0.5** | **0.6** | **0.7** | **0.8** | **0.9** |
|---|---|---|---|---|---|---|---|---|---|---|
| 200 | 0.1 | 0.148 | 0.154 | 0.164 | 0.183 | 0.204 | 0.233 | 0.273 | 0.338 | 0.230 |
| 200 | 0.05 | 0.088 | 0.089 | 0.108 | 0.120 | 0.136 | 0.158 | 0.204 | 0.265 | 0.152 |
| 200 | 0.01 | 0.023 | 0.024 | 0.032 | 0.040 | 0.051 | 0.069 | 0.094 | 0.143 | 0.072 |
| 500 | 0.1 | 0.116 | 0.121 | 0.125 | 0.154 | 0.152 | 0.166 | 0.205 | 0.235 | 0.326 |
| 500 | 0.05 | 0.061 | 0.066 | 0.067 | 0.078 | 0.087 | 0.102 | 0.133 | 0.157 | 0.244 |
| 500 | 0.01 | 0.013 | 0.016 | 0.017 | 0.019 | 0.024 | 0.036 | 0.051 | 0.058 | 0.012 |
| 1000 | 0.1 | 0.120 | 0.100 | 0.110 | 0.120 | 0.130 | 0.150 | 0.168 | 0.200 | 0.230 |
| 1000 | 0.05 | 0.072 | 0.048 | 0.052 | 0.058 | 0.062 | 0.070 | 0.096 | 0.116 | 0.152 |
| 1000 | 0.01 | 0.016 | 0.004 | 0.004 | 0.004 | 0.008 | 0.010 | 0.018 | 0.022 | 0.072 |

TABLE 2
*Power of $V_{n,1}$ when $\delta = 0.2$, $\Delta = 0.5$ and $k^\star = n/2$*

| | | ρ | | | | | | | | |
|---|---|---|---|---|---|---|---|---|---|---|
| **n** | **Level** | **0.1** | **0.2** | **0.3** | **0.4** | **0.5** | **0.6** | **0.7** | **0.8** | **0.9** |
| 200 | 0.1 | 0.586 | 0.528 | 0.474 | 0.420 | 0.379 | 0.366 | 0.357 | 0.385 | 0.488 |
| 200 | 0.05 | 0.467 | 0.414 | 0.367 | 0.325 | 0.289 | 0.271 | 0.270 | 0.301 | 0.411 |
| 200 | 0.01 | 0.261 | 0.223 | 0.192 | 0.165 | 0.147 | 0.131 | 0.137 | 0.166 | 0.260 |
| 500 | 0.1 | 0.868 | 0.806 | 0.720 | 0.623 | 0.525 | 0.431 | 0.359 | 0.311 | 0.349 |
| 500 | 0.05 | 0.779 | 0.669 | 0.608 | 0.509 | 0.408 | 0.326 | 0.263 | 0.223 | 0.267 |
| 500 | 0.01 | 0.551 | 0.453 | 0.369 | 0.282 | 0.212 | 0.158 | 0.122 | 0.104 | 0.137 |
| 1000 | 0.1 | 0.974 | 0.940 | 0.892 | 0.836 | 0.700 | 0.558 | 0.424 | 0.358 | 0.308 |
| 1000 | 0.05 | 0.940 | 0.892 | 0.838 | 0.710 | 0.592 | 0.452 | 0.306 | 0.272 | 0.216 |
| 1000 | 0.01 | 0.818 | 0.700 | 0.596 | 0.492 | 0.362 | 0.254 | 0.148 | 0.128 | 0.118 |

TABLE 3
*Power of $V_{n,1}$ when $\delta = 0.2$, $\Delta = 1$ and $k^\star = n/2$*

| | | ρ | | | | | | | | |
|---|---|---|---|---|---|---|---|---|---|---|
| **n** | **Level** | **0.1** | **0.2** | **0.3** | **0.4** | **0.5** | **0.6** | **0.7** | **0.8** | **0.9** |
| 200 | 0.1 | 0.962 | 0.929 | 0.880 | 0.807 | 0.718 | 0.613 | 0.523 | 0.465 | 0.511 |
| 200 | 0.05 | 0.919 | 0.872 | 0.805 | 0.720 | 0.616 | 0.516 | 0.427 | 0.381 | 0.430 |
| 200 | 0.01 | 0.775 | 0.694 | 0.604 | 0.505 | 0.4112 | 0.325 | 0.265 | 0.230 | 0.277 |
| 500 | 0.1 | 0.999 | 0.999 | 0.993 | 0.976 | 0.935 | 0.840 | 0.686 | 0.494 | 0.406 |
| 500 | 0.05 | 0.999 | 0.993 | 0.981 | 0.947 | 0.877 | 0.756 | 0.577 | 0.398 | 0.321 |
| 500 | 0.01 | 0.980 | 0.957 | 0.905 | 0.819 | 0.702 | 0.533 | 0.356 | 0.230 | 0.182 |
| 1000 | 0.1 | 1 | 1 | 1 | 0.998 | 0.992 | 0.954 | 0.846 | 0.592 | 0.350 |
| 1000 | 0.05 | 1 | 1 | 1 | 0.994 | 0.974 | 0.918 | 0.756 | 0.496 | 0.290 |
| 1000 | 0.01 | 1 | 0.996 | 0.988 | 0.966 | 0.890 | 0.776 | 0.536 | 0.302 | 0.144 |

TABLE 4
*Power of $V_{n,1}$ when $\delta = 0.2$, $\Delta = 1.5$ and $k^\star = n/2$*

| | | ρ | | | | | | | | |
|---|---|---|---|---|---|---|---|---|---|---|
| **n** | **Level** | **0.1** | **0.2** | **0.3** | **0.4** | **0.5** | **0.6** | **0.7** | **0.8** | **0.9** |
| 200 | 0.1 | 0.999 | 0.997 | 0.989 | 0.970 | 0.925 | 0.842 | 0.719 | 0.594 | 0.553 |
| 200 | 0.05 | 0.996 | 0.987 | 0.972 | 0.937 | 0.871 | 0.764 | 0.630 | 0.501 | 0.471 |
| 200 | 0.01 | 0.969 | 0.938 | 0.889 | 0.812 | 0.706 | 0.572 | 0.440 | 0.341 | 0.327 |
| 500 | 0.1 | 1 | 1 | 1 | 1 | 0.997 | 0.977 | 0.904 | 0.717 | 0.493 |
| 500 | 0.05 | 1 | 1 | 1 | 0.998 | 0.989 | 0.952 | 0.838 | 0.622 | 0.402 |
| 500 | 0.01 | 1 | 0.999 | 0.995 | 0.982 | 0.941 | 0.837 | 0.661 | 0.410 | 0.248 |
| 1000 | 0.1 | 1 | 1 | 1 | 1 | 1 | 0.998 | 0.984 | 0.860 | 0.572 |
| 1000 | 0.05 | 1 | 1 | 1 | 1 | 1 | 0.998 | 0.964 | 0.788 | 0.460 |
| 1000 | 0.01 | 1 | 1 | 1 | 1 | 0.996 | 0.970 | 0.876 | 0.606 | 0.270 |

TABLE 5
*Power of $V_{n,1}$ when $\delta = 0.2$, $\Delta = 0.5$ and $k^\star = n/4$*

| | | ρ | | | | | | | | |
|---|---|---|---|---|---|---|---|---|---|---|
| **n** | **Level** | **0.1** | **0.2** | **0.3** | **0.4** | **0.5** | **0.6** | **0.7** | **0.8** | **0.9** |
| 200 | 0.1 | 0.563 | 0.512 | 0.461 | 0.415 | 0.378 | 0.350 | 0.346 | 0.387 | 0.492 |
| 200 | 0.05 | 0.443 | 0.396 | 0.349 | 0.308 | 0.275 | 0.256 | 0.256 | 0.294 | 0.414 |
| 200 | 0.01 | 0.229 | 0.198 | 0.174 | 0.149 | 0.130 | 0.121 | 0.124 | 0.159 | 0.268 |
| 500 | 0.1 | 0.836 | 0.772 | 0.687 | 0.591 | 0.499 | 0.405 | 0.341 | 0.299 | 0.347 |
| 500 | 0.05 | 0.744 | 0.662 | 0.567 | 0.474 | 0.383 | 0.308 | 0.246 | 0.216 | 0.260 |
| 500 | 0.01 | 0.496 | 0.411 | 0.331 | 0.258 | 0.195 | 0.139 | 0.108 | 0.097 | 0.138 |
| 1000 | 0.1 | 0.964 | 0.946 | 0.894 | 0.812 | 0.698 | 0.556 | 0.412 | 0.300 | 0.268 |
| 1000 | 0.05 | 0.942 | 0.890 | 0.814 | 0.706 | 0.580 | 0.428 | 0.310 | 0.226 | 0.176 |
| 1000 | 0.01 | 0.798 | 0.694 | 0.592 | 0.450 | 0.328 | 0.228 | 0.150 | 0.090 | 0.082 |

Table 6
*Power of $V_{n,1}$ when $\delta = 0.2$, $\Delta = 1$ and $k^\star = n/4$*

| | | $\rho$ | | | | | | | | |
|---|---|---|---|---|---|---|---|---|---|---|
| $n$ | **Level** | **0.1** | **0.2** | **0.3** | **0.4** | **0.5** | **0.6** | **0.7** | **0.8** | **0.9** |
| 200 | 0.1 | 0.948 | 0.912 | 0.857 | 0.785 | 0.694 | 0.591 | 0.504 | 0.463 | 0.517 |
| 200 | 0.05 | 0.901 | 0.846 | 0.781 | 0.693 | 0.589 | 0.494 | 0.415 | 0.375 | 0.436 |
| 200 | 0.01 | 0.739 | 0.659 | 0.565 | 0.474 | 0.385 | 0.308 | 0.250 | 0.228 | 0.291 |
| 500 | 0.1 | 0.999 | 0.997 | 0.987 | 0.968 | 0.911 | 0.808 | 0.657 | 0.491 | 0.404 |
| 500 | 0.05 | 0.997 | 0.988 | 0.972 | 0.927 | 0.846 | 0.721 | 0.550 | 0.391 | 0.318 |
| 500 | 0.01 | 0.968 | 0.934 | 0.873 | 0.781 | 0.653 | 0.495 | 0.344 | 0.223 | 0.182 |
| 1000 | 0.1 | 1 | 1 | 1 | 0.998 | 0.994 | 0.952 | 0.820 | 0.590 | 0.368 |
| 1000 | 0.05 | 1 | 1 | 0.998 | 0.996 | 0.980 | 0.904 | 0.730 | 0.464 | 0.278 |
| 1000 | 0.01 | 0.998 | 0.994 | 0.962 | 0.874 | 0.726 | 0.488 | 0.266 | 0.126 | 0.082 |

Table 7
*Power of $V_{n,1}$ when $\delta = 0.2$, $\Delta = 1.5$ and $k^\star = n/4$*

| | | $\rho$ | | | | | | | | |
|---|---|---|---|---|---|---|---|---|---|---|
| $n$ | **Level** | **0.1** | **0.2** | **0.3** | **0.4** | **0.5** | **0.6** | **0.7** | **0.8** | **0.9** |
| 200 | 0.1 | 0.999 | 0.995 | 0.984 | 0.962 | 0.909 | 0.830 | 0.712 | 0.598 | 0.566 |
| 200 | 0.05 | 0.993 | 0.983 | 0.968 | 0.922 | 0.854 | 0.753 | 0.624 | 0.505 | 0.483 |
| 200 | 0.01 | 0.961 | 0.923 | 0.871 | 0.795 | 0.682 | 0.556 | 0.429 | 0.340 | 0.343 |
| 500 | 0.1 | 1 | 1 | 1 | 0.999 | 0.994 | 0.966 | 0.874 | 0.695 | 0.499 |
| 500 | 0.05 | 1 | 1 | 0.999 | 0.995 | 0.982 | 0.931 | 0.808 | 0.601 | 0.402 |
| 500 | 0.01 | 0.999 | 0.997 | 0.990 | 0.966 | 0.915 | 0.803 | 0.622 | 0.409 | 0.239 |
| 1000 | 0.1 | 1 | 1 | 1 | 1 | 1 | 0.998 | 0.974 | 0.826 | 0.472 |
| 200 | 0.05 | 1 | 1 | 1 | 1 | 1 | 0.994 | 0.944 | 0.732 | 0.386 |
| 200 | 0.01 | 1 | 1 | 1 | 1 | 0.992 | 0.952 | 0.814 | 0.522 | 0.230 |

Table 8
*Power of $V_{n,1}$ when $\omega = 1$, $\alpha = 0.1$ $\beta = 0.1$ and $k^* = n/2$*

| $n$ | **Level** | $\Delta$ | **Sim. sign. lev.** | $\Delta$ | **Power** | $\Delta$ | **Power** | $\Delta$ | **Power** |
|---|---|---|---|---|---|---|---|---|---|
| 200 | 0.1 | 0 | 0.182 | 0.5 | 0.615 | 1 | 0.968 | 1.5 | 0.999 |
| 200 | 0.05 | 0 | 0.07 | 0.5 | 0.499 | 1 | 0.962 | 1.5 | 0.997 |
| 200 | 0.01 | 0 | 0.02 | 0.5 | 0.270 | 1 | 0.802 | 1.5 | 0.975 |
| 500 | 0.1 | 0 | 0.114 | 0.5 | 1 | 1 | 1 | 1.5 | 1 |
| 500 | 0.05 | 0 | 0.06 | 0.5 | 1 | 1 | 1 | 1.5 | 1 |
| 500 | 0.01 | 0 | 0.01 | 0.5 | 0.600 | 1 | 0.984 | 1.5 | 0.999 |
| 1000 | 0.1 | 0 | 0.116 | 0.5 | 0.982 | 1 | 1 | 1.5 | 1 |
| 1000 | 0.05 | 0 | 0.06 | 0.5 | 0.962 | 1 | 1 | 1.5 | 1 |
| 1000 | 0.01 | 0 | 0.01 | 0.5 | 0.850 | 1 | 1 | 1.5 | 1 |

The simulations are based again on the assumption that the $\delta_i$'s are independent standard normal variables. Comparing Tables 8 and 9, we can conclude that the ratio test is working well even for small sample sizes when the size of the change is $\Delta = 0$ (no change), $\Delta = 0.5, 1, 1.5$. The values $\omega = 1, \alpha = 0.1$ and $\beta = 0.1$ correspond to very weak dependence between the observations while the choice $\omega = 0.5, \alpha = 0.1$ and $\beta = 0.7$ corresponds to stronger dependence. In both cases the power is high and the same power was obtained for $k^* = n/2$ and $k^* = n/4$.

## 3. Proofs

*Proof of Theorem 1.1.* We can assume without loss of generality that $\mu = 0$. Let

$$Z_{n,1}(t) = n^{-1/2} \sum_{1 \le i \le nt} \epsilon_i \qquad \text{and} \qquad Z_{n,2}(t) = n^{-1/2} \sum_{nt < i \le n} \epsilon_i$$

TABLE 9
*Power of $V_{n,1}$, when $\omega = 0.5, \alpha = 0.1, \beta = 0.7$ and $k^* = n/2$*

| $n$ | Level | $\Delta$ | Sim. sign. lev. | $\Delta$ | Power | $\Delta$ | Power | $\Delta$ | Power |
|---|---|---|---|---|---|---|---|---|---|
| 200 | 0.1 | 0 | 0.128 | 0.5 | 0.478 | 1 | 0.967 | 1.5 | 0.995 |
| 200 | 0.05 | 0 | 0.071 | 0.5 | 0.372 | 1 | 0.842 | 1.5 | 0.985 |
| 200 | 0.01 | 0 | 0.018 | 0.5 | 0.179 | 1 | 0.637 | 1.5 | 0.916 |
| 500 | 0.1 | 0 | 0.107 | 0.5 | 0.775 | 1 | 0.998 | 1.5 | 1 |
| 500 | 0.05 | 0 | 0.060 | 0.5 | 0.663 | 1 | 0.990 | 1.5 | 1 |
| 500 | 0.01 | 0 | 0.012 | 0.5 | 0.423 | 1 | 0.938 | 1.5 | 0.998 |
| 1000 | 0.1 | 0 | 0.116 | 0.5 | 0.954 | 1 | 1 | 1.5 | 1 |
| 1000 | 0.05 | 0 | 0.064 | 0.5 | 0.892 | 1 | 1 | 1.5 | 1 |
| 1000 | 0.01 | 0 | 0.008 | 0.5 | 0.744 | 1 | 0.998 | 1.5 | 1 |

Condition (1.2) yields that

$$(Z_{n,1},(t)Z_{n,2}(t)) \overset{\mathcal{D}^2[0,1]}{\longrightarrow} \sigma(W(t),W^*(t)), \tag{3.1}$$

where $W^*(t) = W(1) - W(t)$. Since

$$n^{-1/2} \sup_{1<i\le nt} \left| \sum_{1\le j\le i} (X_j - \bar{X}_{[nt]}) \right| = \sup_{0<i\le nt} \left| Z_{n,1}(i/n) - \frac{i}{[nt]} Z_{n,1}(t) \right|$$

and similarly

$$n^{-1/2} \sup_{nt<i\le n} \left| \sum_{i\le j\le n} (X_j - \tilde{X}_{[nt]}) \right| = \sup_{nt<i\le n} \left| Z_{n,2}(i/n) - \frac{n-1}{n-[nt]} Z_{n,2}(t) \right|.$$

by (3.1) we have for all $0 < \delta < 1/2$,

$$\left( n^{-1/2} \sup_{0<i\le nt} \left| \sum_{1\le j\le i} (X_j - \bar{X}_{[nt]}) \right|, n^{-1/2} \sup_{nt<i\le n} \left| \sum_{i\le j\le n} (X_j - \tilde{X}_{[nt]}) \right| \right)$$
$$\overset{\mathcal{D}^2[\delta,1-\delta]}{\longrightarrow} \sigma \left( \sup_{0<s\le t} |W(s) - (s/t)W(t)|, \sup_{t\le s\le 1} |W^*(s) - ((1-s)/(1-t))W^*(t)| \right).$$

Hence the proof of (1.3) is complete. The statistics $V_{n,2}$ and $V_{n,3}$ are also continuous functionals of $Z_{n,1}(t), Z_{n,2}(t), 0 \le t \le 1$. Hence the arguments in the proof of (1.3) can be repeated. □

*Proof of Theorem 1.2.* Let $k > k^*$. Then the definition of $X_j$ gives

$$\sum_{j=1}^{i}(X_n - \bar{X}_k) = \begin{cases} \sum_{j=1}^{i} \epsilon_j - \frac{i}{k} \sum_{j=1}^{k} \epsilon_j - \frac{i(k-k^*)}{k}\Delta_n, & \text{if } 1 \le i \le k^* \\ \sum_{j=1}^{i} \epsilon_j - \frac{i}{k} \sum_{j=1}^{k} \epsilon_j + (i-k^*)\Delta_n - \frac{i(k-k^*)}{k}\Delta_n, & \text{if } k^* < k. \end{cases}$$

If $k = [n\tau]$ with some $\theta < \tau < 1 - \delta$, we get that

$$\begin{aligned} n^{-1/2} \max_{1\le i\le k} \left| \sum_{1\le j\le i} (X_j - \overline{X}_k) \right| &\ge n^{-1/2} \left| \sum_{1\le j\le k^*} (X_j - \bar{X}_k) \right| \\ &= O_P(1) + \frac{k^*(k-k^*)}{k} n^{-1/2} |\Delta_n|. \end{aligned}$$

Since there is no change in the means of $X_k, X_{k+1} \dots, X_n$, by Theorem 1.1 we have that

$$n^{-1/2} \max_{k \le i \le n} \left| \sum_{i \le j \le n} (X_j - \tilde{X}_k) \right| \xrightarrow{\mathcal{D}} \sigma(1-\tau)^{1/2} \sup_{0 \le t \le 1} |B(t)|,$$

where $\{B(t), 0 \le t \le 1\}$ is a Brownian bridge. Observing that $k^*(k-k^*)n^{-1/2}|\Delta_n|/k \to \infty$ we conclude that $V_{n,1} \xrightarrow{P} \infty$. Similar arguments yield the proof when $i = 2$ and 3. □

*Proof of Theorem 1.3.* It follows from condition (1.2) that

$$\left\{ n^{-1/2} \max_{1 \le i \le k^*} \left| \sum_{1 \le j \le i} (X_j - \bar{X}_{k^*} \right|, n^{-3/2} \max_{k^* < i \le n} \left| \sum_{i \le j \le n} (X_j - \tilde{X}_k) \right| \right\}$$

$$\xrightarrow{\mathcal{D}} \sigma \left\{ \sup_{0 \le t \le \theta} \left| W(t) - \frac{t}{\theta} W(\theta) \right|, \right.$$
$$\left. \sup_{\theta \le t \le 1} \left| \int_t^1 (W(s) - W(\theta)) ds - \frac{1-t}{1-\theta} \int_\theta^1 (W(s) - W(\theta)) ds \right| \right\},$$

proving that $Z_{n,1} \xrightarrow{P} \infty$. Similar arguments give that $Z_{n,i} \xrightarrow{P} \infty$ when $i = 2$ and 3. □

*Proof of Theorem 1.4.* Simple modifications of of the proof of Theorem 1.3 gives the results. □

*Proof of Theorem 1.5.* It is an immediate consequence of Theorems 1.3 and 1.4. □

## References

[1] Berkes, I., Horváth, L., Kokoszka, P. and Shao, Q.-M. (2005). Almost sure convergence of the Bartlett estimator. *Period. Math. Hungar.* **51** 11–25. MR2180630

[2] Berkes, I., Horváth, L., Kokoszka, P. and Shao, Q.-M. (2006). On discriminating between long–range dependence and changes in the mean. *Ann. Statist.* **34** 1140–1165. MR2278354

[3] Csörgő, M. and Horváth, L. (1997). *Limit Theorems in Change-Point Analysis*. Wiley, Chichester.

[4] Csörgő, M. and Révész, P. (1981). *Strong Approximations in Probability and Statistics*. Academic Press, New York. MR0666546

[5] Giraitis, L., Kokoszka, P., Leipus, R. and Teyssière, G. (2003). Rescaled variance and related tests for long memory in volatility and levels. *J. Econometrics* **112** 265–294. MR1951145

[6] Hannan, E. J. (1979). The central limit theorem for time series regression. *Stochastic Process. Appl.* **9** 281–289. MR0562049

[7] Kim, J.-Y. (2000). Detection of change in persistence of a linear time series. *J. Econometrics* **95** 97–116. MR1746619

[8] Kim, J.-Y., Belaire-Franch, J. and Amador, R. (2002). Corrigendum to "Detection of change in persistence of a linear time series" **109** 389–392.

[9] Leybourne, S. and Taylor, A. (2006). Persistence change tests and shifting stable autoregressions. *Economics Letters* **91** 44–49.

[10] Lo, A. (1991). Long-term memory in stock market prices. *Econometrica* **59** 1279–1313.

[11] Perron, P. (2006). Dealing with structural breaks. In *Palgrave Handbook of Econometrics* **1**. *Econometric Theory* 278–352. Palgrave Mcmillan.

[12] Wang, Q., Lin, Y.-X. and Gulati, C. M. (2002). The invariance principle for linear processes with applications. *Econometric Theory* **18** 119–139. MR1885353

[13] Wu, W. B. and Min, W. (2005). On linear processes with dependent innovations. *Stochastic Process. Appl.* **115** 939–958. MR2138809